\documentclass[12pt,reqno]{amsart}
\usepackage{amsmath}
\usepackage{amsthm}
\usepackage{eurosym}
\usepackage{amssymb}
\usepackage{amscd}
\usepackage{hyperref}
\usepackage[usenames]{color}
\usepackage{mathabx}

\newcommand{\func}{\operatorname}

\theoremstyle{plain}
\newtheorem{theorem}{Theorem}[section]

\theoremstyle{definition}

\newtheorem{remark}[theorem]{Remark}

\newtheorem{example}[theorem]{Example}

\numberwithin{equation}{section}
\begin{document}
\title{Some remarks on equivariant elliptic operators and their invariants}
\author[J.~Br\"{u}ning]{Jochen Br\"{u}ning}
\address{Institut f\"{u}r Mathematik \\
Humboldt Universit\"{a}t zu Berlin \\
Unter den Linden 6 \\
D-10099 Berlin, Germany}
\email[J.~Br\"{u}ning]{bruening@mathematik.hu-berlin.de}
\author[K.~Richardson]{Ken Richardson}
\address{Department of Mathematics \\
Texas Christian University \\
Fort Worth, Texas 76129, USA}
\email[K.~Richardson]{k.richardson@tcu.edu}
\subjclass{58J20; 58J28; 58J35; 57S15}
\keywords{equivariant, index, eta invariant}
\thanks{Work of the first author was partly supported by the grant SFB647
''Space-Time-Matter''. Work of the second author was partly supported by a
grant from the Simons Foundation (Grant Number 245818 to Ken Richardson).}
\dedicatory{Dedicated to the memory of Boris Sternin and Franz Kamber}

\date{\today }
\maketitle

\section{Introduction}

In this note we consider first order elliptic differential operators acting
on smooth vector bundles over compact manifolds, and certain invariants
derived from the analysis of these operators, namely the \textbf{eta
invariant} (in the case of self-adjoint operators) and the \textbf{index}.
While these topics are very well known and established (see \cite{APS1}, 
\cite{APS2}), this is not so for their equivariant counterparts (as
described in \cite{A}). We will develop an approach that works for both the
eta invariant and the index, and we will apply it to the eta invariant in
Section {\ref{equivariantEtaSection}, to the index in Section \ref%
{equivariantIndexSection}, and to the Atiyah-Patodi-Singer Theorem in
Section \ref{equivAPSThmSection}. So far, index formulas in geometric terms
are not available in the general case, but we will give some simple examples
in Section \ref{examplesSection}. }

We note that other researchers have considered the equivariant eta invariant
before. Specifically, H. Donnelly in \cite{Don} provided a fixed point
formula for the equivariant index on a manifold with boundary; this formula
included an equivariant version of the eta invariant dependent on a single
element of the group. In further work, S. Goette \cite{Goette} provided an
infinitesimal version of this theorem (where the formulas depend on an
element of the Lie algebra) and gave applications to equivariant eta
forms in fiber bundles. We mention also the seminal exposition on the
equivariant index in \cite{Be-G-V} and associated work \cite{Be-V1}, \cite%
{Be-V2}. 

What makes this work different from the above is that we are evaluating
integer-valued indices corresponding to multiplicities of group
representations, and our eta invariant is a number dependent on the entire
group at once. In simple cases, such as when the operator is elliptic and
when the group is finite, it is clear how to determine one version of the
equivariant index from the other (and one version of the equivariant eta
invariant from the other). However, for general compact Lie groups and
transversally elliptic operators, our formulas are distinct from the others
obtained before. Moreover, the techniques of proof and formulas obtained are
entirely different and depend on equivariant heat asymptotics that may
involve logarithmic terms.

This note is expository, giving references to all relevant sources. For
simplicity, we consider only elliptic differential operators, even though
the proofs outlined apply to the more general situation. In every case, we
outline the well-known proofs and theorems without Lie group actions first
and then show how these same ideas can be applied in the equivariant cases
with appropriate modifications. A more detailed and expanded article that
applies to transversally elliptic operators will appear in due time.

\section{The equivariant eta invariant\label{equivariantEtaSection}}

Let $\left( M,h_{M}\right) $ be a compact Riemannian manifold of dimension $%
m $, $\left( E,h_{E}\right) $ a Hermitian vector bundle over $M$, and let $%
D_{c}$ be a first order, symmetric elliptic differential operator acting on $%
C^{\infty }\left( M,E\right) $. Since 
\begin{equation*}
\left\langle D_{c}s_{1},s_{2}\right\rangle _{L^{2}\left( M,E\right)
}=\left\langle s_{1},D_{c}s_{2}\right\rangle _{L^{2}\left( M,E\right) }
\end{equation*}%
for $s_{i}\in C^{\infty }\left( M,E\right) $, $i=1,2$, $D_{c}$ is closable;
we write $D:=\overline{D_{c}}$. Then $\func{dom}D=H^{1}\left( M,E\right) $,
the Sobolev space of order $1$, and the fundamental elliptic estimate \cite[%
Thm. III.5.2]{LM} shows that the norm of $H^{1}\left( M,E\right) $, $%
\left\Vert s\right\Vert _{1}$, and the norm $\left\Vert s\right\Vert
_{0}+\left\Vert Ds\right\Vert _{0}$ are equivalent. By Rellich's Theorem 
\cite[Thm. III.2.6]{LM}, we deduce further that $\func{spec}D$ consists of a
discrete set, and any $\lambda \in \func{spec}D$ is an eigenvalue of finite
multiplicity, with the only accumulation points being $\pm \infty $.
Finally, all eigenfunctions are in $C^{\infty }\left( M,E\right) $, by
elliptic regularity \cite[Thm. III.4.5]{LM}.

We now restrict attention to operators with the property that the symbol of $%
D^{2}$ takes the form%
\begin{equation}
\sigma \left( D^{2}\right) \left( x,\xi \right) =\left\vert \xi \right\vert
^{2}I_{E}.  \label{scalarSymbolConditionEq}
\end{equation}%
Also, $D$ is a Fredholm operator, by the parametrix construction. Then%
\begin{equation}
C^{\infty }\left( M,E\right) =\ker D_{c}\oplus D_{c}\left( C^{\infty }\left(
M,E\right) \right) .  \label{smoothDecomp}
\end{equation}%
({\ref{smoothDecomp}}) implies%
\begin{equation*}
L^{2}\left( M,E\right) =\ker D\oplus \func{im}D,
\end{equation*}%
since $\ker D=\ker D_{c}$ and $\func{im}D$ is closed. Since $D\subset
D^{\ast }$ by symmetry, 
\begin{equation*}
L^{2}\left( M,E\right) =\ker D\oplus \func{im}D^{\ast },
\end{equation*}%
and so $\func{im}D=\func{im}D^{\ast }$; hence $D=D^{\ast }$ is self-adjoint.

According to Weyl's Law, we obtain the asymptotics of the eigenvalue
counting function as%
\begin{equation}
N_{D^{2}}\left( t\right) \underset{t\rightarrow \infty }{\sim }C\left(
M,E\right) t^{m/2}.  \label{CountingFunction}
\end{equation}

The invariant we are trying to understand is the \textbf{eta function}, $%
\eta _{D}$, defined as 
\begin{equation*}
\eta _{D}\left( z\right) =\sum_{\lambda \in \mathrm{\func{spec}}D\setminus
\left\{ 0\right\} }\func{sgn}\left( \lambda \right) \left\vert \lambda
\right\vert ^{-z},\quad \operatorname{Re} z\gg 0.
\end{equation*}%
We relate this to $D$ by the Mellin transform%
\begin{equation}
\Gamma \left( z+\frac{1}{2}\right) \eta _{D}\left( 2z\right)
=\int_{0}^{\infty }\func{Tr}\nolimits_{L^{2}}\left[ De^{-tD^{2}}\right]
t^{\left( z+\frac{1}{2}\right) -1}dt,  \label{etaMeromorphExtension}
\end{equation}%
where 
\begin{equation*}
\func{Tr}\nolimits_{L^{2}}\left[ De^{-tD^{2}}\right] =\sum_{\lambda \in 
\mathrm{\func{spec}}D}\lambda e^{-t\lambda ^{2}}=:\alpha _{D}\left( t\right)
.
\end{equation*}%
Obviously, $\alpha _{D}\left( t\right) \underset{t\rightarrow \infty }{%
\longrightarrow }0$ exponentially such that the integral over $[1,\infty )$
is holomorphic in $z\in \mathbb{C}$. A simple calculation using ({\ref%
{CountingFunction}}) shows that $\alpha _{D}\left( t\right) =\mathcal{O}%
\left( t^{-(m+1)/2}\right) $ as $t\rightarrow 0^{+}$ such that $\eta _{D}$
is holomorphic for $\func{Re}z>\frac{m+1}{2}$. To reveal the meromorphic
nature of $\eta _{D}$ in $\mathbb{C}$, we have to study the kernel of $%
De^{-tD^{2}}$ using 
\begin{equation*}
\alpha _{D}\left( t\right) =\int_{M}\func{Tr}\nolimits_{E}\left[
De^{-tD^{2}}\left( x,x\right) \right] ~dx,
\end{equation*}%
and the Hadamard expansion of the kernel of $e^{-tD^{2}}$ (c.f. \cite[Ch. 8,
Thm. 8.1]{Dui}). Then we obtain for $y$ near $x\in M$, trivializing the
vector bundle $E$, 
\begin{equation}
e^{-tD^{2}}\left( x,y\right) \underset{t\rightarrow 0^{+}}{\sim }\left( 4\pi
t\right) ^{-m/2}e^{-d\left( x,y\right) ^{2}\diagup 4t}\sum_{i=0}^{\infty
}t^{i/2}U_{i}\left( x,y\right) .  \label{heatExpansion}
\end{equation}%
Here, $d$ is the Riemannian distance function, and the $U_{i}$ are smooth
sections of $E\boxtimes E^{\ast }$ over $M\times M,$ carrying substantial
geometric information. It is not difficult to check with ({\ref%
{heatExpansion}}), integrating only over $[0,1]$ in ({\ref%
{etaMeromorphExtension}}), that $\eta _{D}$ is meromorphic with simple poles
at the points $z_{i}:=\frac{m-i+1}{2}$.

While the material so far is not new, the situation is different if we look
at the equivariant case. We assume now that a compact Lie group $G$ acts on $%
\left( M,h_{M}\right) $ and $\left( E,h_{E}\right) $ effectively, smoothly,
and isometrically. In addition, we assume that the operator $D_{c}$ commutes
with the $G$-action on $C^{\infty }\left( M,E\right) $. Then we pick an
irreducible representation $\rho $ of $G$, such that the $\rho $-isotypical
subspace $L^{2}\left( M,E\right) _{\rho }$ has infinite dimension. The
orthogonal projection $P_{\rho }$ onto $L^{2}\left( M,E\right) _{\rho }$ is
constructed via the invariant integral, to wit%
\begin{equation}
P_{\rho }s=d_{\rho }\int_{G}\overline{\chi }_{\rho }\left( g\right) ~g\cdot
s~dg=:\int_{G}\kappa _{\rho }\left( g\right) ~g\cdot s~dg~;
\label{projectionFormulaDef}
\end{equation}%
here $d_{\rho }$ is the dimension of a representative space of $\rho $, and $%
\chi _{\rho }$ its character. By assumption, $P_{\rho }$ commutes with $D$
such that 
\begin{equation*}
D_{\rho }=P_{\rho }D=DP_{\rho }=P_{\rho }DP_{\rho }
\end{equation*}%
is a self-adjoint operator, with spectral resolution $\left( P_{\rho
}E_{D}\left( \lambda \right) \right) _{\lambda \in \mathbb{R}}$, since the
spectral resolution of $D$ can be constructed via the resolvent; c.f. \cite[%
VI, Lemma 5.6]{Kato}.

Therefore, $\func{spec}D_{\rho }$ is a subset of $\func{spec}D$, and the
same is true for the corresponding eigenfunctions. Hence we can generalize
formula ({\ref{etaMeromorphExtension}}) as

\begin{eqnarray}
\Gamma \left( z+\frac{1}{2}\right) \eta _{D_{\rho }}\left( 2z\right)
&=&\int_{0}^{\infty }\func{Tr}\nolimits_{L_{\rho }^{2}}\left[ D_{\rho
}e^{-tD_{\rho }^{2}}\right] t^{\left( z+\frac{1}{2}\right) -1}dt  \notag \\
&=&\int_{0}^{\infty }\sum_{\lambda \in \mathrm{\func{spec}}D_{\rho }}\lambda
e^{-t\lambda ^{2}}\cdot t^{\left( z+\frac{1}{2}\right) -1}dt.
\label{etaFunctionInTermsOfKernel_eq}
\end{eqnarray}%
We put%
\begin{equation*}
\alpha _{D_{\rho }}\left( t\right) :=\func{Tr}\nolimits_{L_{\rho }^{2}}\left[
D_{\rho }e^{-tD_{\rho }^{2}}\right]
\end{equation*}%
and obtain the same exponential vanishing as $t\rightarrow \infty $. The
asymptotics for $t\rightarrow 0^{+}$ are more involved; they can be found in 
\cite[Cor. 3.2]{BrH1}. In fact, if $m_{G}:=\dim M_{0}\diagup G$ and $M_{0}$
is the space of principal orbits, then we have%
\begin{equation*}
\alpha _{D_{\rho }}\left( t\right) =\mathcal{O}\left( t^{-(m_{G}+1)/2}\right)
,~t\rightarrow 0^{+}.
\end{equation*}%
The finer analysis concerning the kernel of $D_{\rho }e^{-tD_{\rho }^{2}}$
is essentially contained in \cite[Thm. 4]{BrH2}; we give an explanation of
the main steps.

We choose a cut-off function $\psi \in C^{\infty }\left( \mathbb{R}_{\geq
0},[0,1]\right) $ with $\psi \left( x\right) =1$ for $0\leq x\leq
\varepsilon ^{2}$ and $\psi \left( x\right) =0$ for $x\geq 4\varepsilon ^{2}$%
, and we put $\widetilde{\psi }\left( x,y\right) :=\psi \left(
d^{2}(x,y)\right) $. Then we fix $y$ and let $x\in B_{2\varepsilon }\left(
y\right) $. We also fix a local orthonormal frame $\left( s_{j}\right)
_{j=1}^{\func{rk}E}$ for $\left. E\right\vert _{B_{2\varepsilon }\left(
y\right) }$. Then we obtain, with $u_{ijk}\in C^{\infty }\left(
B_{2\varepsilon }\left( y\right) \right) $,%
\begin{eqnarray}
\widetilde{\psi }\left( x,y\right) e^{-tD^{2}}\left( x,y\right) &=&%
\widetilde{\psi }\left( x,y\right) e^{-d^{2}\left( x,y\right) \diagup 4t}{%
\Biggm [}\sum_{\substack{ -m_G-1\leq i\leq L  \\ 1\leq j,k\leq \func{rk}E}}%
t^{i/2}u_{ijk}\left( x,y\right) s_{j}\left( x\right) \otimes s^{k}\left(
y\right)  \notag \\
&&+\mathcal{O}\left( t^{\frac{L+1}{2}}\right) {\Biggm ]}.  \label{cutoff}
\end{eqnarray}%
We find (using ({\ref{projectionFormulaDef}}) )%
\begin{eqnarray*}
\widetilde{\psi }\left( P_{\rho }e^{-tD^{2}}\right) \left( x,y\right) &=&%
\widetilde{\psi }\left( e^{-tD_{\rho }^{2}}\right) \left( x,y\right) \\
&=&\int_{G}\widetilde{\psi }\left( g\cdot x,y\right) \kappa _{\rho }\left(
g\right) e^{-d^{2}\left( g\cdot x,y\right) \diagup 4t}\cdot \\
&&\cdot \left[ \sum_{i\le L,j,k}t^{i/2}u_{ijk}\left( g\cdot x,y\right)
s_{j}\left( g\cdot x\right) \otimes s^{k}\left( y\right) +\mathcal{O}\left(
t^{\frac{L+1}{2}}\right) \right] ~dg
\end{eqnarray*}%
and 
\begin{eqnarray*}
\func{Tr}\nolimits_{E}\left[ \widetilde{\psi }e^{-tD_{\rho }^{2}}\left(
y,y\right) \right] &=&\int_{G}\widetilde{\psi }\left( g\cdot y,y\right)
\kappa _{\rho }\left( g\right) e^{-d^{2}\left( g\cdot y,y\right) \diagup
4t}\cdot \\
&&\cdot \left[ \sum_{i\le L,j}t^{i/2}u_{ijj}\left( g\cdot y,y\right) +\mathcal{O}%
\left( t^{\frac{L+1}{2}}\right) \right] ~dg.
\end{eqnarray*}%
This gives, finally, with $\widetilde{u}_{ijk}\left( g,g\cdot y,y\right) =%
\widetilde{\psi }\left( g\cdot y,y\right) \kappa _{\rho }\left( g\right)
u_{ijj}\left( g\cdot y,y\right) $,%
\begin{eqnarray*}
&&\func{Tr}\nolimits_{E}\left[ \widetilde{\psi }e^{-tD_{\rho }^{2}}\right] \\
&=&\int_{G\times M}e^{-d^{2}\left( g\cdot y,y\right) \diagup 4t}\cdot \left[
\sum_{i,j}t^{i/2}\widetilde{u}_{ijj}\left( g,g\cdot y,y\right) +\mathcal{O}%
\left( t^{\frac{L+1}{2}}\right) \right] ~~dg~dy.
\end{eqnarray*}%
The crucial part comes from the set%
\begin{equation*}
\mathcal{L}=\{(g,y)\in G\times M:g\cdot y=y\},
\end{equation*}%
since outside any open neighborhood of $\mathcal{L}$, contributions to the
trace are $\mathcal{O}(t^{\infty })$. We choose such a neighborhood of $%
\mathcal{L}$ where $\widetilde{\psi }=1$. It is shown in \cite{BrH2} that
one can construct a \textquotedblleft weak resolution\textquotedblright\ of $%
d^{2}\left( g\cdot y,y\right)$ near $\mathcal{L}$ by a
polynomial in $m$ real variables, the exponential asymptotics of which are readily
computed. We quote the result:

\begin{theorem}
\cite[Thm. 4]{BrH2}\label{equivariantTraceExpansionThm} We have the
asymptotic expansion%
\begin{equation}
\func{Tr}\nolimits_{L^{2}}\left[ e^{-tD_{\rho }^{2}}\right] =\sum_{\substack{
i=-m_{G}  \\ 0\leq j\leq \mathcal{T}(M,G)}}^{L}t^{i/2}\left( \log t\right)
^{j}a_{\rho }^{ij}+\mathcal{O}\left( t^{\frac{L+1}{2}}\right) .
\label{equivariantTraceExpansion}
\end{equation}
\end{theorem}

\begin{remark}
Given the complicated structure of $\mathcal{L}$, the expansion ({\ref%
{equivariantTraceExpansion}}) is essentially qualitative. The number $%
\mathcal{T}(M,G)$ can be taken as 
\begin{equation*}
\mathcal{T}(M,G)=\#\left\{ \dim G\cdot x:x\in M\right\} -1.
\end{equation*}%
As far as we know, no example with non-vanishing $\log $-terms has been
constructed yet. From \cite[Thm. 4]{BrH2} we do know only that $a_{\rho
}^{ij}=0$ for $i=-m_{G},~j>0$.
\end{remark}

\begin{remark}
To apply the asymptotic expansion to the eta invariant calculation, we would
have to add a factor $D_{\rho }$ in ({\ref{equivariantTraceExpansion}}). But
any operator commuting with $e^{-tD_{\rho }^{2}}$ would only change the
coefficients in the expansion and hence will not change the
qualitative structure of the result. In our more detailed work, we show that
for a new set of constants $b_{\rho }^{ij}$,%
\begin{equation*}
\func{Tr}\nolimits_{L^{2}}\left[ D_{\rho }e^{-tD_{\rho }^{2}}\right] =\sum 
_{\substack{ i=-m_{G}-1  \\ 0\leq j\leq \mathcal{T}(M,G)}}^{L}t^{i/2}\left(
\log t\right) ^{j}b_{\rho }^{ij}+\mathcal{O}\left( t^{\frac{L+1}{2}}\right)
\end{equation*}
\end{remark}

\begin{remark}
\label{SaturatedSetExpansion}Using the same calculations and techniques, we
have similar asymptotic expansions of integrals of equivariant heat kernels
over $G\times U$, where $U$ is an open saturated ($G$-invariant) set in $M$; see 
\cite{BrIntegrals}, \cite{BrH2}.
Specifically, 
\begin{equation*}
\int_{G\times U}\widetilde{\psi }\left( g\cdot y,y\right) \kappa _{\rho
}(g)~e^{-tD^{2}}\left( g\cdot y,y\right) ~dg~dy=\sum_{\substack{ i=-m_{G} 
\\ 0\leq j\leq \mathcal{T}(M,G)}}^{L}t^{i/2}\left( \log t\right)
^{j}a_{U,\rho }^{ij}+\mathcal{O}\left( t^{\frac{L+1}{2}}\right)
\end{equation*}
\end{remark}

The result we were after is

\begin{theorem}
$\eta _{D_{\rho }}$ is meromorphic in $\mathbb{C}$ with poles at most in the points $\{-%
\frac{i+1}{2}:i\geq -m_G-1\}$ that may be multiple, also at $0$.
\end{theorem}

\section{The equivariant index\label{equivariantIndexSection}}

\vspace{1pt}Consider $\left( M,h_{M}\right) $ as in Section {\ref%
{equivariantEtaSection}, and introduce two smooth vector bundles }$E^{\pm }$
of rank $\ell $, with Hermitian metrics $h_{E}^{\pm }$. We consider first
order elliptic differential operators%
\begin{equation*}
D_{c}^{\pm }:C^{\infty }\left( M,E^{\pm }\right) \rightarrow C^{\infty
}\left( M,E^{\mp }\right) .
\end{equation*}%
We construct a smooth vector bundle $E=E^{+}\oplus E^{-}$ with metric $%
h_{E}=h_{E}^{+}\oplus h_{E}^{-}$ and a first order symmetric operator%
\begin{equation*}
D_{c}=\left( 
\begin{array}{cc}
0 & D_{c}^{-} \\ 
D_{c}^{+} & 0%
\end{array}%
\right) :\left( 
\begin{array}{c}
C^{\infty }\left( M,E^{+}\right) \\ 
C^{\infty }\left( M,E^{-}\right)%
\end{array}%
\right) \rightarrow \left( 
\begin{array}{c}
C^{\infty }\left( M,E^{+}\right) \\ 
C^{\infty }\left( M,E^{-}\right)%
\end{array}%
\right) .
\end{equation*}%
Then we obtain a self-adjoint Fredholm operator, $D$, as in Section {\ref%
{equivariantEtaSection}, acting on }$H^{1}\left( M,E\right) $, such that
with $D^{\pm }$ the closures of $D_{c}^{\pm }$%
\begin{equation*}
D=\left( 
\begin{array}{cc}
0 & D^{-} \\ 
D^{+} & 0%
\end{array}%
\right) .
\end{equation*}%
We denote by $P^{\pm }$ the orthogonal projections from $E$ to $E^{\pm }$,
and we introduce the self-adjoint involution%
\begin{equation*}
\alpha :=P^{+}-P^{-}\in C^{\infty }\left( M,E\otimes E^{\ast }\right) .
\end{equation*}%
Then we obtain the well-known McKean-Singer formula:%
\begin{eqnarray*}
\func{ind}D^{+} &=&\func{Tr}\nolimits_{L^{2}\left( M,E\right) }\left[ \alpha
e^{-tD^{2}}\right] \\
&=&\func{Tr}\nolimits_{L^{2}\left( M,E\right) }\left[
e^{-tD^{-}D^{+}}-e^{-tD^{+}D^{-}}\right] .
\end{eqnarray*}%
As before, assume that%
\begin{equation*}
\sigma \left( D^{2}\right) \left( x,\xi \right) =|\mathcal{\xi }|^{2}.
\end{equation*}%
We can bring in the group action as in the previous chapter and obtain 
\begin{equation*}
\func{ind}D_{\rho }^{+}=\func{Tr}\nolimits_{L^{2}\left( M,E\right) _{\rho }}%
\left[ \alpha _{\rho }e^{-tD_{\rho }^{2}}\right] ,
\end{equation*}%
So we have to proceed as in the proof of Theorem {\ref%
{equivariantTraceExpansionThm}}. Then we obtain

\begin{theorem}
\begin{equation*}
\func{ind}D_{\rho }=a_{\rho }^{00,+}-a_{\rho }^{00,-},
\end{equation*}%
where the coefficients are constructed in analogy with the procedure leading
to Theorem {\ref{equivariantTraceExpansionThm}}.
\end{theorem}

\section{The Equivariant APS Theorem\label{equivAPSThmSection}}

We now consider a compact manifold with boundary, 
\begin{equation*}
\overline{M}=M\sqcup N,
\end{equation*}%
where $\left( M,h_{M}\right) $ is an open Riemannian manifold with compact
boundary $\left( N,h_{N}\right) $ and $h_{\overline{M}}=h_{M}\sqcup h_{N}$.
We consider also a smooth Hermitian vector bundle $\left(
E_{N},h_{E_{N}}\right) $ over $N$, and a first order symmetric elliptic
operator $D_{N,c}$ acting on $C^{\infty }\left( N,E_{N}\right) $. We
assume that $D_{N}^{2}$ satisfies the condition ({\ref%
{scalarSymbolConditionEq}}).

As in Section {\ref{equivariantEtaSection}, we see that }$D_{N}:=\overline{%
D_{N,c}}$ has domain $H^{1}\left( N,E_{N}\right) $ and is self-adjoint with
a purely discrete spectrum. As in ({\ref{CountingFunction}}), we see that
the eigenvalue counting function satisfies%
\begin{equation*}
N_{D_{N}^{2}}\left( t\right) \underset{t\rightarrow \infty }{\sim }C\left(
N,E_{N}\right) t^{n/2},~n=m-1=\dim N.
\end{equation*}%
We fix an orthonormal eigenbasis $\left( \phi _{\lambda }\right) _{\lambda
\in \func{spec}D_{N}}$ , and we find as in Section {\ref%
{equivariantEtaSection},}%
\begin{eqnarray*}
\alpha _{D_{N}}\left( t\right) :&=&\func{Tr}\nolimits_{L^{2}\left(
N,E_{N}\right) }\left[ D_{N}e^{-tD_{N}^{2}}\right] \\
&=&\int_{N}\func{Tr}\nolimits_{E_{N}}\left[ D_{N}e^{-tD_{N}^{2}}\left(
x,x\right) \right] ~dx \\
&=&\int_{N}\func{Tr}\nolimits_{E_{N}}\left[ \sum_{\lambda \in \func{spec}%
D_{N}}\lambda e^{-t\lambda ^{2}}\phi _{\lambda }\left( x\right) \otimes \phi
_{\lambda }\left( x\right) ^{\ast }\right] ~dx.
\end{eqnarray*}%
This implies as before that 
\begin{eqnarray*}
\alpha _{D_{N}}\left( t\right) &\underset{t\rightarrow \infty }{%
\longrightarrow }&0\text{ exponentially,} \\
\alpha _{D_{N}}\left( t\right) &=&\mathcal{O}\left( t^{-(n+1)/2}\right) \text{
as }t\rightarrow 0.
\end{eqnarray*}

\vspace{1pt}Following again the work in Section {\ref{equivariantEtaSection}
and Section \ref{equivariantIndexSection}, we introduce a compact Lie group, 
}$G$,{\ with properties described there, and we pick an irreducible
representation }$\rho $ of $G$ such that $L^{2}\left( N,E_{N}\right) _{\rho
} $ has infinite dimension. We assume again that $D_{N}$ commutes with $G$,
hence with the orthogonal projection $P_{\rho }$ from ({\ref%
{projectionFormulaDef}}); then $D_{N,\rho }=P_{\rho }D_{N}=D_{N}P_{\rho
}=P_{\rho }D_{N}P_{\rho }$, and 
\begin{eqnarray*}
\alpha _{D_{N,\rho }}(t) :&=&\func{Tr}\nolimits_{L_{\rho }^{2}}\left[
D_{N,\rho }e^{-tD_{N,\rho }^{2}}\right] \\
&=&\sum_{\lambda \in \func{spec}D_{N,\rho }}\lambda e^{-t\lambda ^{2}}
\end{eqnarray*}%
with%
\begin{eqnarray*}
\alpha _{D_{N,\rho }}(t)&\underset{t\rightarrow \infty }{\longrightarrow }&0%
\text{ exponentially, and }\\
\alpha _{D_{N,\rho }}(t)&=&\mathcal{O}\left( t^{-\left( n_{G}+1\right)
/2}\right) \text{ as }t\rightarrow 0^{+},
\end{eqnarray*}%
where $n_{G}=\dim N_0\slash G$.

In the next step, we choose an open neighborhood $U_{N}$ of $N$ in $%
\overline{M}$ such that 
\begin{equation*}
U_{N}\simeq N\times \lbrack 0,\varepsilon )
\end{equation*}%
with product metric $h_{U_{N}}\simeq h_{N}\oplus du^{2}$. We extend $N\times
\lbrack 0,\varepsilon )$ to the cylinder $N\times \mathbb{R}_{\geq 0}=:C_{N}$%
, and we lift the vector bundle $(E_{N},h_{E_{N}})$ to $%
(E_{C_{N}},h_{C_{N}}) $. On $C_{N}$, we introduce the first order elliptic
differential operators 
\begin{equation}
D_{C_{N}}:=\frac{\partial }{\partial u}+D_{N}~,~D_{C_{N}}^{\dag }:=-\frac{%
\partial }{\partial u}+D_{N}~,  \label{D_C_N and Adjoint}
\end{equation}%
both acting on $C^{\infty }\left( C_{N},E_{C_{N}}\right) $. As shown in \cite%
[(2.9), (2.10)]{APS1}, one obtains the identity%
\begin{equation}
\left\Vert s\right\Vert _{1}^{2}=C_{1}\left\Vert s_{0}\right\Vert
_{0}^{2}+C_{2}\left\Vert D_{C_{N}}^{(\dag )}s\right\Vert _{0}^{2}~,~s\in
C^{\infty }\left( C_{N},E_{C_{N}}\right) .  \label{SobolevEquivNorm}
\end{equation}%
Decomposing $s\in C^{\infty }\left( C_{N},E_{C_{N}}\right) $ as%
\begin{equation}
s\left( y,u\right) =\sum_{\lambda \in \func{spec}D_{N}}s_{\lambda }\left(
u\right) \phi _{\lambda }\left( y\right) ,  \label{decompOverCylinder}
\end{equation}%
such that $D_{C_{N}}^{(\dag )}s(y,u)=\sum_{\lambda \in \func{spec}%
D_{N}}\left( (-)s_{\lambda }^{\prime }\left( u\right) +\lambda s_{\lambda
}\left( u\right) \right) \phi _{\lambda }\left( y\right) $, and introducing
the boundary conditions for $D_{C_{N}}$ and $D_{C_{N}}^{\dag }$:%
\begin{eqnarray}
s_{\lambda }\left( 0\right) &=&0~\text{for }\lambda \geq 0~\left(
D_{C_{N}}\right) ,  \label{firstBC} \\
s_{\lambda }\left( 0\right) &=&0~\text{for }\lambda <0~\left(
D_{C_{N}}^{\dag }\right) ,  \label{adjBC}
\end{eqnarray}%
we see that%
\begin{equation*}
\ker D_{C_{N}}=\ker D_{C_{N}}^{\dag }=0.
\end{equation*}%
Abusively, we write $D_{C_{N}}:=\overline{D_{C_{N}}}$ and find from ({\ref%
{SobolevEquivNorm}}) that%
\begin{equation*}
\func{dom}D_{C_{N}}=H^{1}\left( C_{N},E_{N}\right) ,~\ker D_{C_{N}}=0.
\end{equation*}%
As shown in \cite[Prop. 2.42]{APS1}, we have $\overline{D_{C_{N}}^{\dag }}%
=D_{C_{N}}^{\ast }$ and $\func{dom}D_{C_{N}}^{\ast }=H^{1}\left(
C_{N},E_{N}\right) $, too. By the elliptic parametrix construction, both $%
D_{C_{N}}$ and $D_{C_{N}}^{\ast }$ are Fredholm operators with%
\begin{equation}
L^{2}\left( C_{N},E_{C_{N}}\right) =\func{im}D_{C_{N}}=\func{im}%
D_{C_{N}}^{\ast }~,  \label{imagesAreL^2}
\end{equation}%
so both $D_{C_{N}}$ and $D_{C_{N}}^{\ast }$ are isomorphisms.

The operators $D_{C_{N}}$ and $D_{C_{N}}^{\ast }$ generate two Laplacians, $%
\Delta _{C_{N}}^{\pm }$, for which the expression%
\begin{equation}
K\left( t\right) :=\func{Tr}\nolimits_{L^{2}\left( C_{N},E_{N}\right) }\left[
e^{-t\Delta _{C_{N}}^{+}}-e^{-t\Delta _{C_{N}}^{-}}\right]
\label{K(t)_inCylinder}
\end{equation}%
is similar to an index formula but must be of a different nature since $%
C_{N} $ is non-compact. It leads in fact to the surprisingly simple formula (%
\cite[(2.23)]{APS1})%
\begin{equation}
K\left( t\right) =-\sum_{\lambda \in \func{spec}D_{N}}\frac{\func{sgn}%
\lambda }{2}\func{erfc}\left( \left\vert \lambda \right\vert \sqrt{t}\right)
,  \label{K(t)_definition}
\end{equation}%
where the complementary error function $\func{erfc}$ is given by%
\begin{equation}
\func{erfc}(x)=\frac{2}{\sqrt{\pi }}\int_{x}^{\infty }e^{-\xi ^{2}}d\xi <%
\frac{2}{\sqrt{\pi }}e^{-x^{2}};  \label{erfc_stuff}
\end{equation}%
thus, the derivative is given by%
\begin{equation*}
K^{\prime }\left( t\right) =\frac{1}{\sqrt{4\pi t}}\sum_{\lambda \in \func{%
spec}D_{N}}\lambda e^{-\lambda ^{2}t}.
\end{equation*}%
The asymptotic behavior of $K\left( t\right) $ for $t\rightarrow \infty $ is
clear from ({\ref{K(t)_definition}}): if $h:=\dim \ker D_{N}$, then%
\begin{equation*}
K\left( t\right) +\frac{h}{2}\underset{t\rightarrow \infty }{\longrightarrow 
}0\text{ exponentially.}
\end{equation*}%
On the other hand, from ({\ref{K(t)_definition}}) and ({\ref{erfc_stuff}})
and Weyl's Law for $D_{N}$, we find%
\begin{equation*}
\left\vert K\left( t\right) \right\vert \leq \frac{1}{\sqrt{\pi }}%
\sum_{\lambda }e^{-t\lambda ^{2}}\leq Ct^{-n/2}~\text{as }t\rightarrow 0^{+}.
\end{equation*}%
In analogy to the trace formula ({\ref{etaFunctionInTermsOfKernel_eq}}) we
obtain \cite[(2.25)]{APS1}:%
\begin{equation}
-\frac{\Gamma \left( z+\frac{1}{2}\right) }{2z\sqrt{\pi }}\eta
_{D_{N}}\left( 2z\right) =\int_{0}^{\infty }\left( K\left( t\right) +\frac{h%
}{2}\right) t^{z-1}dt.  \label{newEraInTermsOfK(t)_eq}
\end{equation}%
Finally, a careful estimate shows that 
\begin{equation}
K_{U_{N}}\left( t\right) :=\int_{U_{N}}\func{Tr}%
\nolimits_{E_{C_{N}}}[K(t,y,u)]~du~dy=K\left( t\right) +R\left( t\right) ,
\label{approximation by finite tube}
\end{equation}%
where $R\left( t\right) $ is exponentially small as $t\rightarrow 0^{+}$.

Now, we introduce an index problem that is essentially an extension of
the data we have established on $U_{N}$, a collar of the boundary $N$.
Recall that $(\overline{M},h_{\overline{M}})$ is compact with \emph{flat}
collar $U_{N}$. Then we double the manifold to the compact manifold $(%
\widetilde{M},h_{\widetilde{M}})$, such that $(\widetilde{M},h_{\widetilde{M}%
})\supset (\overline{M},h_{\overline{M}})$. Next we introduce two smooth
vector bundles $E$ and $F$ over $\overline{M}$, and a first order elliptic
differential operator%
\begin{equation*}
D:C^{\infty }\left( \overline{M},E\right) \rightarrow C^{\infty }\left( 
\overline{M},F\right) .
\end{equation*}%
We assume that in $U_{N}$,%
\begin{equation}
D=\sigma \left( \frac{\partial }{\partial u}+D_{N}\right) ;
\label{form near boundary}
\end{equation}%
here, $\sigma \in \mathcal{U}\left( E,F\right) $ is the symbol of $D$, $%
\sigma =\sigma \left( du\right) $, and an isometry. We also impose the
boundary condition ({\ref{firstBC}}). Thus, we can construct a (right)
parametrix for $D$ as in Section {\ref{equivariantEtaSection}, and find that 
}$\overline{D}=:D$ (abusively) is a Fredholm operator, and $D^{\ast }$ as
well, with the adjoint boundary condition ({\ref{adjBC}}).

We turn to the index calculation,%
\begin{equation*}
\func{ind}D=\dim \ker D-\dim \ker D^{\ast },
\end{equation*}%
where $\func{dom}D=H^{1}\left( \overline{M},E,P\right) $ and $\func{dom}%
D^{\ast }=H^{1}\left( \overline{M},F,I-P\right) $ such that the operators $%
D^{\ast }D$ and $DD^{\ast }$ have purely discrete spectrum, using Rellich
again. Therefore, the McKean-Singer argument gives again%
\begin{equation*}
\func{ind}D=\func{Tr}\nolimits_{L^{2}\left( \overline{M},E\right) }\left[
e^{-tD^{\ast }D}\right] -\func{Tr}\nolimits_{L^{2}\left( \overline{M}%
,F\right) }\left[ e^{-tDD^{\ast }}\right] .
\end{equation*}%
In the collar $U_{N}$ we can use the function $K(t)$ from ({\ref%
{K(t)_definition}}) as before since the asymptotics of the full integral and
its truncation to $U_{N}$ are the same. Thus we choose a cut-off function $%
\psi $ that equals $1$ near $N$ and has its support in $U_{N}$. We now
extend $D$ to $\widetilde{D}$ acting on $C^{\infty }\left( \widetilde{M},%
\widetilde{E}\right) $ where $\widetilde{E}$ extends $E$, likewise for $%
\widetilde{D}^{\ast }$ and $\widetilde{F}$, and obtain a kernel%
\begin{equation*}
F\left( t\right) =\int_{\widetilde{M}}F\left( t,x\right) ~dx=\func{Tr}%
\nolimits_{L^{2}\left( \widetilde{M},\widetilde{E}\right) }\left[ e^{-t%
\widetilde{D}^{\ast }\widetilde{D}}\right] -\func{Tr}\nolimits_{L^{2}\left( 
\widetilde{M},\widetilde{F}\right) }\left[ e^{-t\widetilde{D}\widetilde{D}%
^{\ast }}\right] .
\end{equation*}%
Now observe that in $U_{N}$, the operators $\widetilde{D}^{\ast }\widetilde{D%
}$ and $\widetilde{D}\widetilde{D}^{\ast }$ are conjugate under $\sigma $
such that their local contribution is zero. Hence $\psi $ can be replaced by 
$1$, and we obtain 
\begin{equation*}
\func{ind}D\underset{t\rightarrow 0}{\thicksim }K\left( t\right) +\int_{%
\overline{M}}F\left( t,x\right) ~dx,
\end{equation*}%
or%
\begin{eqnarray}
&&K\left( t\right) \underset{t\rightarrow 0}{\thicksim }\func{ind}D-\int_{%
\overline{M}}F\left( t,x\right) ~dx  \notag \\
&&\underset{t\rightarrow 0}{\thicksim }\func{ind}D-\sum_{i\geq
-m}t^{i/2}\int_{\overline{M}}c_{i}\left( x\right) ~dx.
\label{K(t)_asympt_formula}
\end{eqnarray}%
Then from ({\ref{newEraInTermsOfK(t)_eq}) and }({\ref{K(t)_asympt_formula})
we obtain}%
\begin{eqnarray}
\int_{0}^{1}\left( K\left( t\right) +\frac{h}{2}\right) t^{z-1}dt &=&-\frac{%
\Gamma \left( z+\frac{1}{2}\right) }{2z\sqrt{\pi }}\eta _{D_{N}}\left(
2z\right) \\
&=&\frac{1}{z}\left( \frac{h}{2}+\func{ind}D\right) -\sum_{i=-m}^{N}\frac{1}{%
\frac{i}{2}+z}\int_{\overline{M}}c_{i}\left( x\right) ~dx+\Theta _{N}(z),
\label{etaInTermsOfKernelFormula}
\end{eqnarray}%
or%
\begin{eqnarray*}
\eta _{D_{N}}\left( 2z\right) &=&-\frac{2z\sqrt{\pi }}{\Gamma \left( z+\frac{%
1}{2}\right) }\left( z^{-1}\left( \frac{h}{2}+\func{ind}D\right)
-\sum_{i=-m}^{N}\frac{c_{i}}{\frac{i}{2}+z}+\Theta _{N}(z)\right) \\
&=&-\frac{2\sqrt{\pi }}{\Gamma \left( z+\frac{1}{2}\right) }\left( \frac{h}{2%
}+\func{ind}D-\sum_{i=-m}^{N}\frac{zc_{i}}{\frac{i}{2}+z}+z\Theta
_{N}(z)\right) ,
\end{eqnarray*}%
or%
\begin{equation*}
-\frac{\Gamma \left( z+\frac{1}{2}\right) }{2\sqrt{\pi }}\eta _{D_{N}}\left(
2z\right) =\frac{h}{2}+\func{ind}D-\sum_{i=-m}^{N}\frac{zc_{i}}{\frac{i}{2}+z%
}+z\Theta _{N}(z).
\end{equation*}%
Here, $c_{i}=\int_{\overline{M}}c_{i}\left( x\right) ~dx$, and $\Theta
_{N}\left( z\right) $ is holomorphic for $\func{Re}z>-\frac{N+1}{2}$. Making 
$z=0$ gives the celebrated APS index formula:%
\begin{equation*}
\func{ind}D=c_{0}-\frac{1}{2}\left( h+\eta _{D_{N}}\left( 0\right) \right) .
\end{equation*}

\vspace{1pt}Finally, we turn to the equivariant case. As before, we
introduce a compact Lie group, $G$, that acts effectively, smoothly and
isometrically on $(N,h_{N})$ and $(E_{N},h_{E_{N}})$. Next we require that
the $G$-action commutes with $D_{N}$ on $C^{\infty }(N,E_{N})$, hence with
the projection $P_{\rho }$ onto $L^{2}(C_{N},E_{N})_{\rho }$. Then 
\begin{equation*}
D_{N,\rho }:=P_{\rho }D_{N}=D_{N}P_{\rho }=P_{\rho }D_{N}P_{\rho }.
\end{equation*}%
As before, we choose $\rho \in \widehat{G}$ such that $\dim
L^{2}(C_{N},E_{N})_{\rho }=\infty $, and observe that for $s\in C^{\infty
}(N,E_{N})$,%
\begin{equation*}
D_{N,\rho }s=\sum_{\lambda \in \func{spec}D_{N,\rho }}\lambda \left\langle
s,\phi _{\lambda }\right\rangle \phi _{\lambda },
\end{equation*}%
such that $\func{spec}D_{N,\rho }\subset \func{spec}D_{N}$. Hence all
equivariant calculations concerning the analogue, $K_{\rho }(t)$, of $K(t)$
are obtained by replacing $\func{spec}D_{N}$ with $\func{spec}D_{N,\rho }$.
To describe this, we extend the $G$-action to $C_{N}$. Since these actions
are isometric, they extend to the cylinder $C_{N}=N\times \mathbb{R}_{\geq
0} $ with trivial action on $\mathbb{R}_{\geq 0}$. The operators ({\ref%
{D_C_N and Adjoint}}) commute with this action, too, with common domain $%
H^{1}(C_{N},E_{C_{N}})_{\rho }$. The decomposition ({\ref{decompOverCylinder}%
}) and the properties of the operators $D_{C_{N},\rho }$ and $D_{C_{N},\rho
}^{\ast }$, including the boundary conditions, are analogous to those
described earlier in this section.

Then we can form the Laplacians%
\begin{equation*}
\Delta _{C_{N},\rho }^{+}:=D_{C_{N},\rho }^{\ast }D_{C_{N},\rho },~~\Delta
_{C_{N},\rho }^{-}:=D_{C_{N},\rho }D_{C_{N},\rho }^{\ast },
\end{equation*}%
and the trace analogous to ({\ref{K(t)_inCylinder}}),%
\begin{equation*}
K_{\rho }(t)=\func{Tr}\nolimits_{L^{2}\left( C_{N},E_{C_{N}}\right) _{\rho }}%
\left[ e^{-t\Delta _{C_{N},\rho }^{+}}-e^{-t\Delta _{C_{N},\rho }^{-}}\right]
,
\end{equation*}%
that is not expressing an index.

Furthermore, we obtain the formula%
\begin{equation*}
-\frac{\Gamma \left( z+\frac{1}{2}\right) }{2z\sqrt{\pi }}\eta _{D_{N},\rho
}\left( 2z\right) =\int_{0}^{1}\left( K_{\rho }\left( t\right) +\frac{%
h_{\rho }}{2}\right) t^{z-1}dt,
\end{equation*}%
where%
\begin{equation*}
h_{\rho }=\dim P_{\rho }\ker D_{N}.
\end{equation*}%
As in ({\ref{approximation by finite tube}}), we obtain 
\begin{equation*}
K_{U_{N},\rho }(t)-K_{C_{N},\rho }(t)\underset{t\rightarrow 0^{+}}{%
\longrightarrow }0\text{ exponentially.}
\end{equation*}

Now we bring in smooth $G$-invariant vector bundles $E$ and $F$ over $%
\overline{M}$, and a first order elliptic differential operator, $D$, that
commutes with $G$. We assume that $D$ has the form ({\ref{form near boundary}%
}) such that $\sigma $ commutes with $G$, too. We further assume that $D$ is
closed with domain $H^{1}(\overline{M},E,P)$, and that $D^{\ast }$ has
domain $H^{1}(\overline{M},F,I-P)$. Thus, $D$ and $D^{\ast }$ are again
Fredholm operators that commute with $G$, by the triviality of the $G$%
-action on $\mathbb{R}_{\geq 0}$. Hence, with $\rho \in \widehat{G}$ as
before, we obtain Hilbert spaces $L^{2}(\overline{M},E,P)_{\rho }$ and $%
L^{2}(\overline{M},F,I-P)_{\rho }$ and operators $D_{\rho }$ and $D_{\rho
}^{\ast }$ with domains $H^{1}(\overline{M},E,P)_{\rho }$ and $H^{1}(%
\overline{M},F,I-P)_{\rho }$, respectively. Then we arrive at the index
formula%
\begin{equation*}
\func{ind}D_{\rho }=\func{Tr}\nolimits_{L^{2}(\overline{M},E,P)_{\rho }}%
\left[ e^{-tD_{\rho }^{\ast }D_{\rho }}\right] -\func{Tr}\nolimits_{L^{2}(%
\overline{M},F,I-P)_{\rho }}\left[ e^{-tD_{\rho }D_{\rho }^{\ast }}\right] .
\end{equation*}%
To evaluate this, we need a good parametrix for the two Laplacians. In the
collar $U_{N}$ (as in ({\ref{imagesAreL^2}}) ), a good parametrix can be
constructed from $K_{\rho }(t)$. In fact, there we have%
\begin{equation*}
D^{\ast }D=\Delta _{C_{N}}^{+},~DD^{\ast }=\sigma \Delta _{C_{N}}^{-}\sigma
^{-1},
\end{equation*}%
such that we obtain from ({\ref{K(t)_inCylinder}})%
\begin{equation*}
K_{\rho }(t)\underset{t\rightarrow 0^{+}}{\thicksim }\func{Tr}%
\nolimits_{L^{2}(U_{N},E_{C_{N}})_{\rho }}\left[ e^{-tD_{C_{N},\rho }^{\ast
}D_{C_{N},\rho }}-e^{-tD_{C_{N},\rho }D_{C_{N},\rho }^{\ast }}\right] .
\end{equation*}%
To obtain a good parametrix for the equivariant Laplacians on $\overline{M}%
\setminus U_{N}$, we proceed as we did earlier in this section. We double $M$
to $\widetilde{M}$, the vector bundles $E$ and $F$ to $\widetilde{E}$ and $%
\widetilde{F}$, and also the metrics $h_{\overline{M}}$, $h_{E}$, and $h_{F}$
to $\widetilde{h}_{\widetilde{M}}$, $\widetilde{h}_{E}$, and $\widetilde{h}%
_{F}$. We also double the operators $D$ and $D^{\ast }$ to $\widetilde{D}$
and $\widetilde{D}^{\ast }$ and obtain $G$-equivariant data in all cases.
Now we can apply Theorem {\ref{equivariantTraceExpansionThm} to }$\widetilde{%
D}_{\rho }^{\ast }\widetilde{D}_{\rho }$ and $\widetilde{D}_{\rho }%
\widetilde{D}_{\rho }^{\ast }$ acting on $C^{\infty }(\widetilde{M},%
\widetilde{E})_{\rho }$ and $C^{\infty }(\widetilde{M},\widetilde{F})_{\rho
} $, respectively, such that%
\begin{eqnarray}
\func{ind}\widetilde{D}_{\rho } &=&\func{Tr}\nolimits_{L^{2}(\widetilde{M},%
\widetilde{E})_{\rho }}\left[ e^{-t\widetilde{D}_{\rho }^{\ast }\widetilde{D}%
_{\rho }}\right] -\func{Tr}\nolimits_{L^{2}(\widetilde{M},\widetilde{F}%
)_{\rho }}\left[ e^{-t\widetilde{D}_{\rho }\widetilde{D}_{\rho }^{\ast }}%
\right]  \notag \\
&&\underset{t\rightarrow 0^{+}}{\thicksim }\sum_{\substack{ i\geq -m_{G}  \\ %
0\leq j\leq \mathcal{T}(\widetilde{M},G)}}t^{i/2}\log ^{j}t(\widetilde{a}%
_{\rho }^{+,ij}-\widetilde{a}_{\rho }^{-,ij}).  \label{indexD_rho^tilde}
\end{eqnarray}%
Using Remark {\ref{SaturatedSetExpansion}, we note that if we integrate the
difference of the equivariant heat kernels over }$\overline{M}$ instead of $%
\widetilde{M}$, we obtain a similar formula (with different constants). That
is,%
\begin{equation*}
\int_{\overline{M}}\left( \func{Tr}\nolimits_{E}\left[ e^{-t\widetilde{D}%
_{\rho }^{\ast }\widetilde{D}_{\rho }}\right] -\func{Tr}\nolimits_{E}\left[
e^{-t\widetilde{D}_{\rho }\widetilde{D}_{\rho }^{\ast }}\right] \right)
\left( x,x\right) ~dx\underset{t\rightarrow 0^{+}}{\thicksim }\sum 
_{\substack{ i\geq -m_{G}  \\ 0\leq j\leq \mathcal{T}(\widetilde{M},G)}}%
t^{i/2}\log ^{j}t(\overline{a}_{\rho }^{+,ij}-\overline{a}_{\rho }^{-,ij}).
\end{equation*}%
To derive the formula for $D_{\rho }$ with domain $H^{1}(\overline{M},E_{%
\overline{M}},P)$ we choose a cut-off function $\psi $ as in ({\ref{cutoff}}%
) and construct the parametrix as%
\begin{eqnarray*}
&&K_{\rho }(t,x)\psi (x)+\left( \func{Tr}\nolimits_{E}\left[ e^{-t\widetilde{%
D}_{\rho }^{\ast }\widetilde{D}_{\rho }}\right] -\func{Tr}\nolimits_{E}\left[
e^{-t\widetilde{D}_{\rho }\widetilde{D}_{\rho }^{\ast }}\right] \right)
(x,x)(1-\psi (x)) \\
&=&:K_{\rho }(t,x)\psi (x)+F_{\rho }(t,x)(1-\psi (x)),~x\in M.
\end{eqnarray*}%
Then it is easy to see that $\psi $ can be replaced by $1$, in view of the
asymptotics of $K_{\rho }$, as detailed in ({\ref{K(t)_definition}}) and ({%
\ref{erfc_stuff}}), and the vanishing of $F_{\rho }(t,x)$ on $U_{N}$ since
the relevant operators are conjugated by $\sigma $ there.

To arrive at our final result, we see Section {\ref{equivariantEtaSection}}
and the references there, and the formula%
\begin{equation*}
\func{ind}D_{\rho }=K_{\rho }(t)+F_{\rho }(t),
\end{equation*}%
where 
\begin{equation*}
F_{\rho }(t)=\sum_{\substack{ i\geq -m_{G}  \\ 0\leq j\leq \mathcal{T}(%
\widetilde{M},G)}}t^{i/2}\log ^{j}t(\overline{a}_{\rho }^{+,ij}-\overline{a}%
_{\rho }^{-,ij}).
\end{equation*}%
Then we can use the calculation as in \cite[p. 56]{APS1} to prove the
following with ({\ref{indexD_rho^tilde}}):

\begin{theorem}
\label{equivAPSTheorem}%
\begin{equation*}
\func{ind}D_{\rho }=\overline{a}_{\rho }^{+,00}-\overline{a}_{\rho }^{-,00}-%
\frac{1}{2}(h_{\rho }+\eta _{D_{N},\rho }),
\end{equation*}%
where $h_\rho =\dim \ker D_{N,\rho }$, and $\eta _{D_{N},\rho }$ is
the constant term in the Laurent expansion of the eta function at $0$.
\end{theorem}

\section{Examples\label{examplesSection}}

\subsection{The de Rham-Hodge operator and the equivariant Euler
characteristic}

It is well known that if $M$ is a Riemannian manifold and $f:M\rightarrow M$
is an isometry that is homotopic to the identity, then the Euler
characteristic of $M$ is the sum of the Euler characteristics of the fixed
point sets of $f$. We generalize this result as follows. We consider the de
Rham operator 
\begin{equation*}
d+d^{\ast }:\Omega ^{\func{even}}\left( M\right) \rightarrow \Omega ^{\func{%
odd}}\left( M\right)
\end{equation*}%
on a $G$-manifold, and the invariant index of this operator is the
equivariant Euler characteristic $\chi ^{G}\left( M\right) $, the Euler
characteristic of the elliptic complex consisting of invariant forms. If $G$
is connected and the Euler characteristic is expressed in terms of its $\rho 
$-components, only the invariant part $\chi ^{G}\left( M\right) =\chi ^{%
\mathbf{1}}\left( M\right) $ appears. This is a consequence of the homotopy
invariance of de Rham-Hodge cohomology. Thus $\chi ^{G}\left( M\right) =\chi
\left( M\right) $ for connected Lie groups $G$. In general the Euler
characteristic is a sum of components%
\begin{equation*}
\chi \left( M\right) =\sum_{\left[ \rho \right] }\chi ^{\rho }\left(
M\right) =\sum_{\left[ \rho \right] }\func{ind}\left( \left( d+d^{\ast
}\right) _{\rho }:\Omega _{\rho }^{\func{even}}\left( M\right) \rightarrow
\Omega _{\rho }^{\func{odd}}\left( M\right) \right) .
\end{equation*}%
Here, $\chi ^{\rho }\left( M\right) $ is the alternating sum of the
dimensions of the $\left[ \rho \right] $-parts of the cohomology groups (or
spaces of harmonic forms). Since the connected component $G_{0}$ of the
identity in $G$ acts trivially on the harmonic forms, the only nontrivial
components $\chi ^{\rho }\left( M\right) $ correspond to representations
induced from unitary representations of the finite group $G\slash G_{0}$.

\begin{example}
Let $M=S^{n}$, let $G=O\left( n\right) $ act on latitude spheres
(principal orbits, diffeomorphic to $S^{n-1}$). Then there are two strata,
with the singular stratum being the two poles. Without using the theorem,
since the only harmonic forms are the constants and multiples of the volume
form, we see that%
\begin{equation*}
\chi ^{\rho }\left( S^{n}\right) =\left\{ 
\begin{array}{ll}
\left( -1\right) ^{n} & \text{if }\rho =\xi \\ 
1 & \text{if }\rho =\mathbf{1}%
\end{array}%
\right. ,
\end{equation*}%
where $\xi $ is the induced one dimensional representation of $O\left(
n\right) $ on the volume forms.
\end{example}

\begin{example}
If instead the group $\mathbb{Z}_{2}$ acts on $S^{n}$ by the antipodal map,
note that%
\begin{equation*}
\chi ^{\rho }\left( S^{n}\right) =\left\{ 
\begin{array}{ll}
1-1=0 & \text{if }\rho =\mathbf{1}\text{ and }n\text{ is odd} \\ 
1 & \text{if }\rho =\mathbf{1}\text{ and }n\text{ is even} \\ 
1 & \text{if }\rho =\xi \text{ and and }n\text{ is even} \\ 
0 & \text{otherwise}%
\end{array}%
\right.
\end{equation*}%
since the antipodal map is orientation preserving in odd dimensions and
orientation reversing in even dimensions.
\end{example}

\begin{example}
Consider the action of $\mathbb{Z}_{4}$ on the torus $T^{2}=\mathbb{R}^{2}%
\slash\mathbb{Z}^{2}$, where the action is generated by a $\frac{\pi }{2}$
rotation. Explicitly, $k\in \mathbb{Z}_{4}$ acts on $\left( 
\begin{array}{c}
y_{1} \\ 
y_{2}%
\end{array}%
\right) $ by 
\begin{equation*}
\phi \left( k\right) \left( 
\begin{array}{c}
y_{1} \\ 
y_{2}%
\end{array}%
\right) =\left( 
\begin{array}{cc}
0 & -1 \\ 
1 & 0%
\end{array}%
\right) ^{k}\left( 
\begin{array}{c}
y_{1} \\ 
y_{2}%
\end{array}%
\right) .
\end{equation*}%
Endow $T^{2}$ with the standard flat metric. The harmonic forms have basis $%
\left\{ 1,dy_{1},dy_{2},dy_{1}\wedge dy_{2}\right\} $. Let $\rho _{j}$ be
the irreducible character defined by $k\in \mathbb{Z}_{4}\mapsto e^{ikj\pi
/2}$. Then the de Rham-Hodge operator $\left( d+d^{\ast }\right) ^{\mathbf{1}%
}$ on $\mathbb{Z}_{4}$-invariant forms has kernel $\left\{
c_{0}+c_{1}dy_{1}\wedge dy_{2}:c_{0},c_{1}\in \mathbb{C}\right\} $. One also
sees that $\ker \left( d+d^{\ast }\right) ^{\rho _{1}}=\func{span}\left\{
idy_{1}+dy_{2}\right\} $, $\ker \left( d+d^{\ast }\right) ^{\rho
_{2}}=\left\{ 0\right\} $, and $\ker \left( d+d^{\ast }\right) ^{\rho _{3}}=%
\func{span}\left\{ -idy_{1}+dy_{2}\right\} $. Then 
\begin{equation*}
\chi ^{\mathbf{1}}\left( T^{2}\right) =2,\chi ^{\rho _{1}}\left(
T^{2}\right) =\chi ^{\rho _{3}}\left( T^{2}\right) =-1,\chi ^{\rho
_{2}}\left( T^{2}\right) =0.
\end{equation*}%
This illustrates the point that it is not possible to use the Atiyah-Singer
integrand on the principal stratum to compute even the invariant index
alone. Indeed, the Atiyah-Singer integrand would be a constant times the
Gauss curvature, which is identically zero. In these cases, the three
singular points $a_{1}=\left( 
\begin{array}{c}
0 \\ 
0%
\end{array}%
\right) ~,a_{2}=\left( 
\begin{array}{c}
0 \\ 
\frac{1}{2}%
\end{array}%
\right) ~,a_{3}=\left( 
\begin{array}{c}
\frac{1}{2} \\ 
\frac{1}{2}%
\end{array}%
\right) $ certainly contribute to the index. The quotient $T^{2}\slash%
\mathbb{Z}_{4}$ is an orbifold homeomorphic to a sphere.
\end{example}

\subsection{The Dolbeault operator on complex projective space}

Consider the action of $g=\left( e^{i\theta _{0}},...,e^{i\theta
_{n}}\right) $ in $T^{n+1}$ on complex projective space $\mathbb{CP}^{n}$
defined by 
\begin{equation*}
g\left[ z_{0},z_{1},...,z_{n}\right] =\left[ e^{i\theta
_{0}}z_{0},...,e^{i\theta _{n}}z_{n}\right] .
\end{equation*}%
The orbit types for this action are as follows. Let $S\subset \left\{
0,1,...,n\right\} $ be a subset, and define 
\begin{equation*}
M_{S}=\left\{ \left[ z_{0},...,z_{n}\right] :z_{k}\neq 0\text{ if and only
if }k\in S\right\} .
\end{equation*}%
Then the isotropy subgroup associated to this orbit type is 
\begin{equation*}
H_{S}=\left\{ \left( e^{i\theta _{0}},...,e^{i\theta _{n}}\right) \in
T^{n+1}:e^{i\theta _{j}}=e^{i\theta _{k}}\text{ whenever }j,k\in S\right\} .
\end{equation*}%
The partial ordering on the isotropy subgroups $H_{S}$ corresponds exactly
to the partial ordering $\subset $ on subsets of $\left\{ 0,1,...,n\right\} $%
. Note that $M_{S}$ is a product of cylinders. That is, if $S=\left\{
k_{0},...,k_{p}\right\} $, then we may set $z_{k_{0}}=1$, and then $%
z_{k_{j}}\in \mathbb{C}- \left\{ 0\right\} $ for $j=1,...,p$, so that $M_{S}$
is topologically and complex-analytically a product of $p$ copies of $%
\mathbb{C}- \left\{ 0\right\} $. Also, note that $G=T^{n+1}$ acts by
rotations in each copy of $\mathbb{C}- \left\{ 0\right\} $. Therefore ${%
G\backslash M_S}$ is a $p$-fold product of open half-lines {$%
S^1\backslash\left( \mathbb{C}- \left\{ 0\right\} \right) $}. There are a
total of $n+1$ fixed points with $H_{S}=T^{n+1}$ that comprise the minimal
orbit type; there are $\left( 
\begin{array}{c}
n+1 \\ 
2%
\end{array}%
\right) $ orbit types $M_{S}\cong \mathbb{C}- \left\{ 0\right\} $ with $%
H_{S}\cong T^{n}$; in general, there are $\left( 
\begin{array}{c}
n+1 \\ 
k+1%
\end{array}%
\right) $ orbit types $M_{S}\cong \left( \mathbb{C}- \left\{ 0\right\}
\right) ^{k}$ with $H_{S}\cong T^{n-k+2}$ for $0\leq k\leq n$. Each orbit
type with $H_{S}\cong T^{n-k+2}$ is a $T^{k-1}$-bundle over a $\left(
k-1\right) $-cell, and the orbits are the fibers of this bundle.

We now consider the operator $\overline{\partial }+\overline{\partial }%
^{\ast }:\Omega ^{j,\text{even}}\rightarrow \Omega ^{j,\text{odd}}$ for
some $j$. Since the Dolbeault cohomology $H^{p,q}$ is one-dimensional only
when $p=q$ and is zero otherwise (generated by powers of the K\"{a}hler
form), the index of $\overline{\partial }+\overline{\partial }^{\ast }$ is $%
\left( -1\right) ^{j}$. Since these cohomology classes are represented by
harmonic forms in $H^{2p}\left( \mathbb{CP}^{n}\right) $, they are invariant
by the (connected) torus group action. Thus, 
\begin{equation*}
\func{ind}\nolimits_{\rho }\left( \overline{\partial }+\overline{\partial }%
^{\ast }:\Omega ^{j,\text{even}}\rightarrow \Omega ^{j,\text{odd}}\right)
=\left\{ 
\begin{array}{ll}
\left( -1\right) ^{j}, & \text{if }\rho =\mathbf{1}\text{;} \\ 
0, & \text{otherwise.}%
\end{array}%
\right.
\end{equation*}

\subsection{The equivariant eta invariant of the boundary signature operator}

Let $B=i^{n}\left( -1\right) ^{p+1}\left( \ast d-d\ast \right) $ on $2p$%
-forms be the boundary signature operator on even degree forms on a manifold
of dimension $2n-1$. In \cite{APS2}, the authors consider a generalization $%
\eta _{\alpha }$ of the eta invariant of this operator, where the forms are
twisted by a representation of the fundamental group. Inside the proof of 
\cite[Prop 2.12]{APS2}, the authors consider the lens space $S^{2n-1}\diagup 
\mathbb{Z}_{m}$ with $S^{2n-1}\subset \mathbb{C}^{n}$, where the generator
of $\mathbb{Z}_{m}$ acts on $\mathbb{C}^{n}$ by multipling the $j^{\text{th}%
} $ copy of $\mathbb{C}$ by $\exp \left( i\theta _{j}\right) $, and the
action is free on the sphere. If we wish to calculate for $\tau $ an
irreducible character of $\mathbb{Z}_{m}$ (say the one that takes the
generator to $\exp \left( 2\pi i\ell /m\right) $), then we may note that%
\begin{equation*}
B_{\ell }:=B_{\tau }\cong \left( B\otimes V_{\tau }^{\ast }\right) ^{G},
\end{equation*}%
so that in the \cite{APS2} notation, the equivariant eta invariant
corresponding to the irreducible character $\tau \neq 1$ is 
\begin{eqnarray*}
\eta _{B_{\ell }}\left( 0\right) &=&\eta _{\tau ^{\ast }}\left(
S^{2n-1}\diagup \mathbb{Z}_{m},0\right) =\rho _{\tau ^{\ast }}\left(
S^{2n-1}\diagup \mathbb{Z}_{m}\right) +\eta _{1}\left( S^{2n-1}\diagup 
\mathbb{Z}_{m},0\right) \\
&=&\frac{1}{m}\sum_{g\neq 1}\sigma _{g}\left( S^{2n-1}\right) \left\{ \chi
_{\tau ^{\ast }}\left( g\right) -1\right\} .
\end{eqnarray*}

From the \cite{APS2} calculation, we get that%
\begin{equation*}
\eta _{B_{\ell }}\left( 0\right) =\frac{1}{m}\sum_{k=1}^{m-1}i^{-n}\left(
\prod\limits_{j=1}^{n}\cot \left( \frac{1}{2}k\theta _{j}\right) \right)
\left( \exp \left( -2\pi i\frac{k\ell }{m}\right) -1\right).
\end{equation*}%
Note that $\left( \exp \left( -2i\theta \right) -1\right) =-2ie^{-i\theta }\left( 
\frac{e^{i\theta }-e^{-i\theta }}{2i}\right) =2i^{-1}e^{-i\theta }\sin
\left( \theta \right) ,$ so our formula may also be written as%
\begin{equation*}
\eta _{B_{\ell }}\left( 0\right) =\frac{2i^{-n-1}}{m}\sum_{k=1}^{m-1}\left(
\prod\limits_{j=1}^{n}\cot \left( \frac{1}{2}k\theta _{j}\right) \right)
\exp \left( -\pi i\frac{k\ell }{m}\right) \sin \left( \frac{k\ell }{m}\pi
\right).
\end{equation*}%
It might be surprising to the reader that this number is always real. To
see this, notice that the quantity has the form%
\begin{equation*}
\eta _{B_{\ell }}\left( 0\right) =\frac{1}{m}\sum_{k=1}^{m-1}\left(
\prod\limits_{j=1}^{n}\frac{1+u_{j}^{k}}{1-u_{j}^{k}}\right) \left(
u_{0}^{k}-1\right) ,
\end{equation*}%
where $u_{0}$ and each $u_{j}$ are $m^{\text{th}}$ roots of unity. If we
replace $k$ with $m-k$ in the sum, we obtain the same result but get the
complex conjugate. Thus the formula simplies further: if $n$ is even, 
\begin{eqnarray*}
\eta _{B_{\ell }}\left( 0\right) &=&\frac{2\left( -1\right) ^{n/2+1}}{m}%
\sum_{k=1}^{m-1}\left( \prod\limits_{j=1}^{n}\cot \left( \frac{1}{2}k\theta
_{j}\right) \right) \sin \left( \frac{k\ell }{m}\pi \right) \sin \left( 
\frac{k\ell }{m}\pi \right) \\
&=&\frac{2\left( -1\right) ^{n/2+1}}{m}\sum_{k=1}^{m-1}\left(
\prod\limits_{j=1}^{n}\cot \left( \frac{1}{2}k\theta _{j}\right) \right)
\sin ^{2}\left( \frac{k\ell }{m}\pi \right) .
\end{eqnarray*}%
If $n=2s+1$ is odd,

\begin{eqnarray*}
\eta _{B_{\ell }}\left( 0\right) &=&\frac{2\left( -1\right) ^{s+1}}{m}%
\sum_{k=1}^{m-1}\left( \prod\limits_{j=1}^{n}\cot \left( \frac{1}{2}k\theta
_{j}\right) \right) \cos \left( \frac{k\ell }{m}\pi \right) \sin \left( 
\frac{k\ell }{m}\pi \right) \\
&=&\frac{\left( -1\right) ^{s+1}}{m}\sum_{k=1}^{m-1}\left(
\prod\limits_{j=1}^{n}\cot \left( \frac{1}{2}k\theta _{j}\right) \right)
\sin \left( \frac{2k\ell }{m}\pi \right) .
\end{eqnarray*}

For instance, for $S^{3}$ with $n=2$, $\theta _{1}=\frac{2\pi }{6}$ and $%
\theta _{2}=\frac{10\pi }{6}$ with $m=6$, 
\begin{eqnarray*}
\eta _{B_{\ell }}\left( 0\right)  &=&\frac{1}{3}\sum_{k=1}^{5}\left( \cot
\left( \frac{k\pi }{6}\right) \cot \left( \frac{5k\pi }{6}\right) \right)
\sin ^{2}\left( \frac{k\ell }{6}\pi \right) \\
&=&\left\{ 
\begin{array}{ll}
0 & \ell =0 \\ 
-\frac{2}{3} & \ell =1,5 \\ 
-\frac{5}{3}~ & \ell =2,4 \\ 
-2 & \ell =3%
\end{array} .
\right. 
\end{eqnarray*}


\begin{thebibliography}{99}
\bibitem{A} M. F. Atiyah, \emph{Elliptic operators and compact groups},
Lecture Notes in Mathematics \textbf{401}, Berlin: Springer-Verlag, 1974.

\bibitem{APS1} M. F. Atiyah, V. K. Patodi, and I. M. Singer, \emph{Spectral
asymmetry and Riemannian geometry. I}, Math. Proc. Camb. Phil. Soc. \textbf{%
77} (1975), 43--69.

\bibitem{APS2} M. F. Atiyah, V. K. Patodi, and I. M. Singer, \emph{Spectral
asymmetry and Riemannian geometry. II}, Math. Proc. Camb. Phil. Soc. \textbf{%
78} (1975), 405--432.

\bibitem{Be-G-V} N. Berline, E. Getzler, and M. Vergne, \emph{Heat Kernels
and Dirac operators}, Grundlehren der mathematischen Wissenschaften \textbf{%
298}, Springer-Verlag, Berlin, 1992.

\bibitem{Be-V1} N. Berline and M. Vergne, \emph{The Chern character of a
transversally elliptic symbol and the equivariant index}, Invent. Math. 
\textbf{124}(1996), no. 1-3, 11-49.

\bibitem{Be-V2} N. Berline and M. Vergne, \emph{L'indice \'{e}quivariant des
op\'{e}rateurs transversalement elliptiques}, Invent. Math. \textbf{124}%
(1996), no. 1-3, 51-101.

\bibitem{BrIntegrals} J. Br\"{u}ning, \emph{On the asymptotic expansion of
some integrals}. Arch. Math. \textbf{42} (1984), no. 3, 253--259.

\bibitem{BrH1} J. Br\"{u}ning and E. Heintze, \emph{Representations of
compact Lie groups and elliptic operators}, Inv. Math. \textbf{50} (1979),
169-203.

\bibitem{BrH2} J. Br\"{u}ning and E. Heintze, \emph{The asymptotic expansion
of Minakshisundaram--Pleijel in the equivariant case,} Duke Math. Jour. 
\textbf{51} (1984), 959-979.

\bibitem{Don} H. Donnelly, \emph{Eta invariants for }$G$\emph{-spaces},
Indiana Univ. Math. J. \textbf{27}(1978), no. 6, 889-918.

\bibitem{Dui} J. J. Duistermaat, \emph{The heat kernel Lefschetz fixed point
formula for the spin-c Dirac operator}, Reprint of the 1996 edition, Modern
Birkh\"{a}user Classics, Birkh\"{a}user/Springer, New York, 2011.

\bibitem{Goette} S. Goette, \emph{Equivariant $\eta$-invariants and $\eta$%
-forms}, J. Reine Angew. Math. \textbf{526} (2000), 181--236.

\bibitem{Kato} T. Kato, \emph{Perturbation theory for linear operators},
Classics in Mathematics, Reprint of the 1980 edition, Springer-Verlag,
Berlin, 1995.

\bibitem{LM} H. B. Lawson and M.-L. Michelsohn, \emph{Spin Geometry},
Princeton University Press, Princeton, 1989.


\end{thebibliography}
\end{document}